\documentclass[11pt]{amsart}
\usepackage[bitstream-charter]{mathdesign}
\usepackage[T1]{fontenc}
\usepackage[latin1]{inputenc}

\usepackage{color,cancel,amsmath}

\numberwithin{equation}{section}
\usepackage{hyperref}
\newcommand{\R}{{\mathbb R}}

\def\R{I \!\!R}

\newcommand{\hsob}{H^1_0(\Omega)}
\newcommand{\dhsob}{H^{-1}(\Omega)}
\newcommand{\elle}[1]{L^{#1}(\Omega)}
\def\w2{W_0^{1,2}(\Omega)}
\def\w{W_0^{1,2}(\Omega)}

\def\io{\int\limits_{\Omega}}
\def\norma#1#2{\|#1\|_{\lower 4pt \hbox{$\scriptstyle #2$}}}

\newcommand{\vomega}{\varOmega}
\newcommand{\into}{{\int_{\vomega}}}

\newcommand{\vfi}{\varphi}

\newcommand\um {u_m}

\def\zm{z_m}

\newcommand{\Jm}{J_m}

\newcommand{\refe}[1]{{(\ref{#1})}}
\newcommand{\dys}{\displaystyle}
\newcommand{\lbl}{\label}
\def\be{\begin{equation}}
\def\ee{\end{equation}}

\newtheorem{teo}{Theorem}[section]
\newtheorem{proposition}[teo]{Proposition}

\newtheorem{remark}[teo]{Remark}

\newtheorem{example}{Example}[section]

\newtheorem{prop}[example]{Proposition}

\def\sqr#1#2{{\vcenter{\vbox{\hrule height.#2pt
\hbox{\vrule
width.#2pt height#1pt \kern#1pt\vrule width.#2pt}
\hrule height.#2pt}}}}
\begin{document}
\title[Asymptotic behavior of positive solutions]{Asymptotic behavior of positive solutions of semilinear elliptic problems with increasing powers}
\author[Lucio Boccardo, Liliane Maia \& Benedetta Pellacci]{Lucio Boccardo, Liliane Maia \& Benedetta Pellacci}
 
\address[L. Boccardo]{Istituto Lombardo and Sapienza Universit\`a di Roma,}
\email{boccardo@uniroma1.it}
\address[L. Maia]{Departamento de Matem\'atica, 
UNB Universidade de Bras\'\i lia, 70910-900 Bras\'\i lia, Brazil}
\email[L. Maia]{lilimaia@unb.br}
\address[B. Pellacci]{Dipartimento di Matematica e Fisica, 
Universit\`a della Campania ``Luigi Vanvitelli'', V.le  A.  Lincoln 5, Caserta, Italy.}
\email[B. Pellacci]{benedetta.pellacci@unicampania.it}
\subjclass[2010]{35B40; 35J66, 35J20}
\keywords{Increasing powers in the lower order term, asymptotical analysis, semilinear elliptic equations}
\thanks{Research partially supported by:   PRIN-2017-JPCAPN Grant: ``Equazioni 
differenziali alle derivate parziali non lineari'',
by project Vain-Hopes within the program VALERE: VAnviteLli pEr la RicErca and by the INdAM-GNAMPA group. 
L.\ Maia was partially supported by FAPDF, CAPES, and CNPq grant 308378/2017 -2.}
\maketitle
\begin{abstract}
\noindent
We prove existence results of  two solutions of the problem 
\[
\begin{cases}
L(u)+u^{m-1}=\lambda u^{p-1}  & \text{ in $\Omega$},
\\
\quad u>0 &\text{ in $\Omega$},
\\
\quad u=0 &  \text{ on $\partial \Omega$},
\end{cases}
\]
where $L(v)=-{\rm div}(M(x)\nabla v)$ is a linear operator,
 $p\in (2,2^{*}]$ and $\lambda$ and $ m$ sufficiently large.
Then their asymptotical limit as $m\to +\infty$ is investigated showing different 
behaviors.
\end{abstract}


\section{Introduction} 
Let $\Omega$ be a bounded domain in $\R^{N}$. 
We  study the asymptotical behavior for $m$ tending to infinity of some positive solutions of the
following semi-linear problem
\be  \label{Prm}
\begin{cases}
L(u)+u^{m-1}=\lambda u^{p-1}  & \text{ in $\Omega$},
\\
\quad u=0 &  \text{ on $\partial \Omega$},
\end{cases}
\ee where $L(v)=-{\rm div}(M(x)\nabla v):\hsob\to \dhsob$ is a
linear ope\-rator in divergence form.  The matrix $M(x)=(m_{ij}(x))$
is symmetric, bounded and positive definite, i.e. there exist
positive constants $0<\alpha<\beta$ such that
 \be\label{alfa} \alpha
|\xi|^2 \leq M(x)\xi \cdot\xi \leq \beta |\xi|^2.
 \ee
The  exponents $p$, $m$  such that
\begin{equation}\lbl{pmz}
2<p\leq 2^*<m, \quad \text{where  }  
2^*=\begin{cases}\frac{2N}{N-2} & \text{for $N>2$,}
\\
2^{*}=+\infty &\text{ for $N=1,2$} .
\end{cases}
\end{equation}
and $\lambda$ is a positive parameter. 

In order to perform our asymptotical analysis, we  first study \eqref{Prm} for $m$ large but fixed and we  prove the following result.
\begin{teo}\lbl{esi}
Assume conditions \eqref{alfa} and \eqref{pmz}. Then,   there exists
$\underline{\lambda}>0$ such that for each
$\lambda>\underline{\lambda}$ there is  $m_0>2^*$ such that for
every  $m\geq m_0$, problem \eqref{Prm} has, at least,  two positive
solutions $\um\not\equiv \zm\in \hsob\cap \elle{m}$. 
\end{teo}
The solutions $u_{m}$ and $z_{m}$ are found as
nontrivial critical points of the  functional 
 \begin{equation}\lbl{defJ}
J_m(v):=\frac12 \into M(x)\nabla v\cdot \nabla v
+\frac{1}m\|v\|_m^m-\frac{\lambda}p\|v^+\|_p^p.
\end{equation}
As $m$ does not satisfy any bound from above (see \eqref{pmz}) 
we  are naturally lead to  study $J_{m}$ defined on $H^{1}_{0}(\Omega)\cap \elle{m}$, and we show that $J_{m}$ has a global minimum point $u_{m}$
at a negative level,  and  another critical point, found applying the Ambrosetti-Rabinowitz Theorem at a positive minimax level, that is $z_{m}$.
Even at this  stage we can see some crucial difference between 
these two solutions: the existence of a nontrivial global minimum point
$u_{m}$ follows assuming the weaker hypothesis on $p$
\begin{equation}\lbl{pmu}
2<p <m; 
\end{equation}
while, in order to have the existence of $z_{m}$ at a positive action level
we need \eqref{pmz}.

Then we  perform the asymptotical analysis of the solutions found for $m\to +
\infty$. In this study a crucial role will be played  by an $L^{\infty}(\Omega)$ a priori bound, which will prevent any blow up of the sequences of
$u_{m}$ and $z_{m}$ even when $p=2^{*}$.
This  marks a strong difference with the well-known explosion phenomenon
for the Lane-Emden problem in bounded domain when $p$ approaches
$2^{*}$ for $N\geq3$, as shown in \cite{atkpel,brepel,han,rey}, and even
with the finite asymptotical behavior of the $L^{\infty}(\Omega)$ norm of 
least energy solution in  dimension 2 as obtained in \cite{adigro}.
Here, we see a completely different scenario which is more related with
the results obtained in \cite{dalors} and in \cite{bocmur} in the case of 
a datum $f(x)$ in the place of the nonlinearity $u^{p-1}$.

The absence  of blowing up family  of solution is  due to the fact the diverging exponent is the one that involves the ``absorbing term''
$u^{m-1}$ on the left hand side instead of the ``reaction one'' $u^{p-1}$ on the right hand side. 
 
As a consequence of the a priori bounds we are  able to show 
the following result
\begin{teo}\label{limintro}
There exist $u,\,z\in {\mathcal K}:=\{v\in H^{1}_0(\Omega):0\leq v(x) \leq 1 \},$ such  that $u,\,z\not\equiv 0$, and 
$u_{m}\rightharpoonup u$ weakly in $\hsob$, $z_{m}\to z$ strongly in $\hsob$
and the convergence is for both $u_{m}$ and $z_{m}$ strongly in
every Lebesgue space. 
Moreover, $u$ and $z$ satisfy the following variational inequalities
\begin{align}\label{ineqintro}
 \int_\Omega M(x)\nabla w    \cdot\nabla (v-w)    
 \geq \lambda \io \,w^{p-1}(  v - w ),\quad \forall\, v\in {\mathcal K,}
\end{align}
In addition,  there exist
$g_{u},\,g_{z}\in L^{\infty}(\Omega)$ , such that
\[
0\leq g_{u},\, g_{z}\leq \lambda \quad
g_{u}(x)[1-u(x)]=0,\; g_{z}(x)[1-z(x)]=0 \;\text{  a. e. in $\Omega$},
\]
and it results
\begin{equation}\label{eqlimintro}
L(u)+g_{u}=\lambda u^{p-1},\qquad L(z)+g_{z}=\lambda z^{p-1}, \quad 
\text{in  $\Omega$}.
\end{equation}
\end{teo}
Theorem \ref{limintro} will be proved as a consequence of Theorems
\ref{limmini}, \ref{limcrit} and it shows that Problem
\eqref{Prm} works as a nonlinear penalization  procedure to obtain solutions
of the variational inequality \eqref{ineqintro} (see \cite{dalors, bocmur}).
The role of the functions $g_{u}$ and $g_{z}$ appearing in 
the limit equation \eqref{eqlimintro}  is well described by the 
equalities $g_{u}(1-u)= 0$ and $g_{z}(1-z)= 0$ almost everywhere.
Indeed, these condition tell us that  $g_{u}$ and $g_{z}$ weight the set in which $u\equiv1$ and $z\equiv 1$, respectively.

In the convergence procedure we can see other difference between
the behavior of $u_{m}$ and $z_{m}$:
First of all, we do not need to prove that $u_{m}$ strongly  converges in
$\hsob$  to a non-trivial limit $u$  and even more we get this limit
for every  $p>2$; on the other hand, we need to show that $z_{m}$ strongly
converges to $z$ in $\hsob$ in order to  show  that $z$ is not trivial. Moreover,
this is obtained by  a comparison with the minimax value  of the 
limit functional 
\[
J_{\infty}:=\frac12 \into M(x)\nabla v\cdot \nabla v-\frac{\lambda}p\into |v|^p,
\]
and the assumption $p\leq 2^{*}$  is crucial  at this point
(see for more  details Remarks \ref{strongconv}-\ref{pcriti}).

However, the main difference between the two asymptotical behaviours
concerns the function $g_{u}$ and $g_{z}$. 
Indeed, while it is quite simple to show that
$g_{u}$ is not trivial thanks to an a priori uniform bound on the level of the global 
minimum point $u_{m}$ of $J_{m}$, the analogous analysis  on $g_{z}$  seems
to be really delicate and strongly depending on the domain and on the exponent 
$p$.
For example, when $\Omega$ is star-shaped and $p=2^{*}$ then 
$g_{z}$ must be not trivial and the same occurs
when $\Omega$ is a ball and $p\to 2^{*}$ from below or $p\to 2$ from above
(see Proposition \ref{pcrit}); 
but whether or not $g_{z}$  is trivial  
for a general $p\in (2,2^{*})$ remains an open question:
we can only prove an abstract result  (see Proposition \ref{gnontrivial})  when the domain
$\Omega$ is a ball 
that actually shows  that $g_{z}$ may be trivial or not 
(see for more details Remark \ref{ending}). 

The paper is organized as follows: in Section \ref{fixedm}  we show the existence 
of $u_{m}$ and $z_{m}$ for fixed $m$ (sufficiently large), in 
Section \ref{bounds} the crucial a-priori bounds and in Section \ref{limit}
we perform the asymptotical analysis.


\section{Existence Results for fixed $m$ }\label{fixedm}

We denote with $\|\cdot\|_r$ the norm in the Lebesgue space
$\elle{r}$ for $1\leq r\leq \infty$ and with $\|\cdot\|$ the usual
norm in the Sobolev space $\hsob$. We will study problem \eqref{Prm}
by variational methods, so that we consider the functional space $
X_m:=\hsob\cap \elle{m}$ endowed with the norm
$\|\cdot\|_{X_m}=\|\cdot\|+\|\cdot\|_m$ and the functional 
$J_m\,:\,X_m\to \R$  defined in \eqref{defJ}.
A solution of Problem \refe{Prm} is a nonzero critical point of
$J_m$, i.e. a function $w\in X_m$, $w(x)>0$, a.e. $x\in\Omega$, such
that the following equation is satisfied for every $v\in X_m$
\begin{equation}\lbl{equat}
\into M(x)\nabla w\cdot\nabla v+\into |w|^{m-2}w\,v-\lambda \into |w|^{p-2}w\,v=0
\end{equation}

 We will prove Theorem \ref{esi} as a consequence of  two
 existence results, the first one is concerned with the existence of the minimal
 solution $u_{m}$
 \begin{teo}\label{esimini}
Assume conditions \eqref{alfa} and \eqref{pmu}.
Then,   there exists
$\underline{\lambda}_{1}>0$  such that for each
$\lambda>\underline{\lambda}_{1}$ there is  $m_0>2^*$ such that for
every  $m\geq m_0$, Problem \eqref{Prm} has a nonnegative, minimal solution
solutions $\um \in \hsob\cap \elle{m},\, u_{m}\not \equiv 0$  with 
\begin{equation}\label{stimamin}
J_m(\um)\leq-\delta<0.
\end{equation}
\end{teo}
\begin{remark}
Let us observe that $\underline{\lambda}_{1}>0$ is explicitly given in \eqref{eq:defA}).
\end{remark}
\begin{proof}
 From H\"older and Young
inequalities we obtain
\begin{align*}
J_m(v) & \geq \frac\alpha
2\|v\|^2+\frac{1}m\|v\|_m^m-\frac{\lambda}p\|v\|_{m}^p
|\Omega|^{1-p/m}
\\
& \geq  \frac\alpha 2\|v\|^2+\frac{1}{2m}\|v\|_m^m
-\frac{m-p}{p^{2}m}|\Omega|\lambda^\frac{m}{m-p} 2^{\frac{p}{m-p}},
\end{align*}
which implies that $J_m$ is coercive in $X_m$. Moreover, if
$v_m\rightharpoonup v$ in $X_m$, then $v_m\to v$ almost everywhere (up to a subsequence),
so that $v_m\to v$ strongly in $\elle p$ thanks to Egorov Theorem; this fact, and the weak
lower semicontinuity of the norm in $X_m$ imply that $J_m$ is weakly
lower semicontinuous. Then, there exists a global minimum point
$\um$ of $J_m$. 
In order to prove \eqref{stimamin},  for any  function $\psi\in \hsob\cap L^{\infty}(\Omega)$, we consider the one variable, positive,
real function $g_{m}(t)$ defined by
\begin{equation}\lbl{AB=}
g_{m}(t)=\frac{a}{t^{p-2}}+ b_m t^{m-p},\qquad a= \frac{\beta}{2} \frac{\| \psi\|^2}{\|\psi\|_p^p} \ , \quad b_m=
\frac{1}{m} \frac{\dys \| \psi\|_m^m}{\dys \|\psi\|_p^p} .
\end{equation}
Note that $m>p>2$ implies that  $\lim\limits_{t\to 0}
g(t)=\lim\limits_{t\to +\infty} g(t)=+\infty$ so that, $g$ attains
its global minimum at the point  $T_m$   given by
\begin{equation}\label{eq:tm}
T_m=\left[ \frac{p-2}{m-p} \frac{a}{b_m}  \right]^{\frac{1}{m-2}}=
\left[\frac{\beta}2\frac{m(p-2)}{m-p}\|\psi\|^{2}\right]^{\frac1{m-2}}
\|\psi\|_{m}^{-\frac{m}{m-2}}.
\end{equation}
Notice that 
\begin{equation}\label{MinJ<0}
J_m(T_m\psi)<0\; \Leftrightarrow\; \lambda>\lambda_{m}(\psi):=p\min_{(0,+\infty)}g_{m}=pg_{m}(T_{m})
\end{equation}
where $g_{m}(T_{m})$ is given by
\begin{eqnarray*}
g(T_m)&=&  \left( \frac{p-2}{m-p} \right)^{\frac{2-p}{m-2}}
 \frac{a^{\frac{m-p}{m-2}}  }{b_m^{\frac{2-p}{m-2}}}
 +
\left( \frac{p-2}{m-p} \right)^{\frac{m-p}{m-2}}
\frac{a^{\frac{m-p}{m-2}}}{b_m^{\frac{2-p}{m-2} }}
  \\
   \\
&=& a^{\frac{m-p}{m-2}} b_m^{\frac{p-2}{m-2}} \left[ \left(
\frac{p-2}{m-p} \right)^{\frac{2-p}{m-2}} + \left( \frac{p-2}{m-p}
\right)^{\frac{m-p}{m-2}} \right]
\\
   \\
&=& a^{\frac{m-p}{m-2}} b_m^{\frac{p-2}{m-2}}  \left(
\frac{p-2}{m-p} \right)^{\frac{2-p}{m-2}}  \frac{m-2}{m-p} .
\end{eqnarray*}
Notice that it results
\begin{equation} \label{lim:Qm}
\lim_{m\to+\infty} \lambda_{m}(\psi)=
p \frac{\beta}{2} \frac{\|\psi\|^2
\|\psi\|_\infty^{p-2}}{\|\psi\|_p^p} =: \Lambda(\psi),
\end{equation}
and, taking into account \eqref{AB=} and  \eqref{eq:tm}, one has
\begin{equation}\label{Tmconv}
\begin{split}
\lim_{m\to+\infty}T_m &=
\lim_{m\to+\infty}\left[\frac1{b_{m}}\right]^{\frac1{m-2}}
=\lim_{m\to+\infty}\left[\frac{1}{m} \frac{\dys \| \psi\|_m^m}{\dys \|\psi\|_p^p}\right]^{-\frac1{m-2}}=\frac1{\|\psi\|_{\infty}} .
\end{split}
\end{equation}
 Moreover, again recalling \eqref{AB=}
\begin{equation}\label{eq:tmpsi}
\begin{split}
\lim_{m\to+\infty}\left[T_{m}\|\psi\|_{m}\right]^{m}
&=\lim_{m\to+\infty}\left[ (p-2)\frac{\beta}{2}\frac{m}{m-p}
\right]^{\frac{m}{m-2}}\|\psi\|^{\frac{2m}{m-2}}
\|\psi\|_{m}^{-\frac{2m}{m-2}}
\\
&=\frac{\beta}2(p-2)\|\psi\|^{2}\|\psi\|_{\infty}^{-2}.
\end{split}\end{equation}
This, joint with \eqref{defJ} and \eqref{Tmconv} yields
\[
\lim_{m\to+\infty} J_{m}(T_{m}\psi)=J_{\infty}(T_{\infty}\psi)=\frac{\|\psi\|_{\infty}^{-2}}2
\into M(x)\nabla \psi\cdot\nabla \psi 
-\lambda
\frac{\|\psi\|_{\infty}^{-p}}p\|\psi\|_{p}^{p} ,
\]
where $J_{\infty}$ is defined by
\begin{equation}\label{def:Jinfty}
J_{\infty}(v)=\frac{1}2
\into M(x)\nabla v\cdot\nabla v
-\frac{\lambda}p\|v\|_{p}^{p} .
\end{equation}
Notice that, due to \eqref{pmu} $J_{\infty}$
may not be finite, but we are computing $J_{\infty}$ only on 
$\psi\in \hsob\cap L^{\infty}(\Omega)$. Moreover, it easy to see that
\begin{equation}\label{disjinftymini}
J_{\infty}(T_{\infty}\psi)<0 \;\Leftrightarrow\; \lambda >\Lambda(\psi).
\end{equation}

Then, $\underline{\lambda}_{1}$ given by
\begin{equation}\label{eq:defA} 
\underline{\lambda}_{1}:=\inf_{\psi\in H^{1}_{0}(\Omega)\cap  L^{\infty}(\Omega)}
\Lambda(\psi)=\inf_{\psi\in H^{1}_{0}(\Omega)\cap  L^{\infty}(\Omega)}
p \frac{\beta}{2} \frac{\|\psi\|^2
\|\psi\|_\infty^{p-2}}{\|\psi\|_p^p} 
\end{equation}
is well defined as $\Lambda(\vfi_{1})$ is finite where $\vfi_{1}$ is the 
positive eigenfunction associated with  the first positive eigenvalue of the 
laplacian operator with homogeneous Dirichlet boundary conditions, $\lambda_{1}^{\rm Dir}$.
Moreover, $\underline{\lambda}_{1}$ is positive as 
\[
\begin{split}
\into |\psi^{p}|=\into |\psi|^{p-2}|\psi|^{2}\leq \|\psi\|_{\infty}^{p-2}\into |\psi|^{2}\leq \lambda_{1}^{\rm Dir}\|\psi\|_{\infty}^{p-2}\|\psi\|^{2}
\end{split}\]
so that $\underline{\lambda_{1}}\geq p\frac{\beta}2\lambda_{1}^{\rm Dir}$.
Thus,  for every $\lambda> \underline{\lambda}_{1} $, we can fix
\begin{equation}\label{eq:defpsi}
\psi_{0}\in \hsob\cap\elle{\infty}\qquad\text{such that}\qquad \Lambda(\psi_{0})<\lambda 
\end{equation}
yielding \eqref{disjinftymini}.
Taking into account \eqref{eq:tmpsi}, we can find  $\sigma>0 $ and sufficiently small and 
 $m_1 >2^*$ such that  for $m\geq  m_1 >2^*$
\begin{equation}\label{eq:psi0}
J_{m}(T_{m}\psi_{0})\leq J_{\infty}(T_{\infty}\psi_{0})+\sigma=J_{\infty}\left(\frac{\psi_{0}}{\|\psi_{0}\|_{\infty}}\right)+\sigma<0
\end{equation}
and the conclusion follows choosing $-\delta=J_{\infty}(\psi_{0}/\|\psi_{0}\|_{\infty}).$

Finally, it is possible to obtain $u_{m}\geq 0$ by  considering  the 
modified functional 
\[
J(v)=\frac12\into M(x)\cdot \nabla v\cdot \nabla v+\frac1m\into |v|^{m}-\frac{\lambda}p\into (v^{+})^{p}.
\]
The same argument of the proof of Theorem \ref{esimini} yields
a nonnegative minimum point $u_{m}$. 

\end{proof}
\begin{remark}
Let us  point out  that  the existence of $u_{m}$ does not require  
$m$ sufficiently large.
On the other hand,  the assumption on $m$ sufficiently large
is crucial to obtain the a priori upper bound on the level of the minimum point $u_{m}$.
\\
Moreover,  concerning the exponent $p$, the existence of the minimum point $u_{m}$ also holds  for  $p$  super-critical (with respect of the Sobolev embeddings).
\end{remark}
The existence of the second nontrivial solution is given by the fol\-lo\-wing 
existence result.
 \begin{teo}\label{esimp}
Assume conditions \eqref{alfa} and \eqref{pmz}. Then,   there exists
$\underline{\lambda}\geq \underline{\lambda}_{1}>0$ (see \eqref{eq:defA}) 
such that for each $\lambda>\underline{\lambda} $ there is  $m_0>\max\{2^*,m_{1}\} $
($m_{1}$ introduced in Theorem
\ref{esimini}), such that for
every  $m\geq m_0$, Problem \eqref{Prm} has a nonnegative, critical point 
 $\zm \in \hsob\cap \elle{m}$  with $ \zm\not\equiv 0$ and $\zm\not\equiv \um$.
\end{teo}
 \begin{proof}
We will obtain the existence of $z_{m}$ by applying  the Ambosetti-Rabinowitz
Theorem \cite{amra}.
First of all, note that
\begin{align*}
J_{m}(v)
&\geq
\frac12\|\nabla v\|_{2}^{2}+\frac1m\|v\|_{m}^{m}-
\lambda\frac{\mathcal S^{p/2}}p|\Omega|^{1-p/2^{*}}\|\nabla v\|_{2}^{p}
\\
&=\frac1m\|v\|_{m}^{m}+\|\nabla v\|_{2}^{2}\left(\frac12-\lambda\frac{\mathcal S^{p/2}}p|\Omega|^{1-p/2^{*}}\|\nabla v\|_{2}^{p-2}\right).
\end{align*}
Now, fix $r_{\lambda}$ such that
\[
r_{\lambda}:=\min\left\{1 ,
\left(\frac{p}{4\lambda{\mathcal S^{p/2}}|\Omega|^{1-p/2^{*}}} \right)^{1/p-2}\right\}
\]
and consider $v\in X_{m}$ such that $\|v\|_{X_{m}}=r_{\lambda}$. In
case $r_{\lambda}\leq 1$ one obtains
(for $m\geq 4$)
\begin{equation}\label{weromin}
\begin{split}
J_{m}(v)
&\geq \frac1m\|v\|_{m}^{m}+\frac14 \|\nabla v\|_{2}^{2}
\geq
\frac1{m}\left(\|v\|_{m}^{m}+\|\nabla v\|_{2}^{m}\right)
\\
&\geq
\frac1{m2^{m-1}}\left(\|v\|_{m}+\|\nabla v\|_{2} \right)^{m}
=\frac1{m2^{m-1}}r_{\lambda}^{m}=:\rho_{m,\lambda}.
\end{split}\end{equation}
In this way we have proved that there exists $r_{\lambda}$ and $\rho_{m,\lambda}$ such that
\[
J_{m}(v)\geq \rho_{m,\lambda}\qquad \forall v\in X_{m}: \|v\|_{X_{m}}=r_{\lambda}.
\]
Then, we can consider the family of paths
$$
\Gamma_m=\left\{\gamma:[0,1]\to X_m,\,:\gamma \text{ is continuous
and } \gamma(0)=0, \,\gamma(1)=T_{m}\psi _{0} \right\}.
$$
where $T_{m}$ and $\psi_{0}$ are defined in \eqref{eq:tm}
and \eqref{eq:defpsi} respectivley.
and
$$
c_m:=\inf_{\Gamma_m}\max_{[0,1]}J_m(\gamma(t)).
$$
Notice that in order to have $T_m\|\psi_{0}\| > r_\lambda $
it is sufficient to have  that
$$
\lambda >
\frac{\alpha}{4\sigma} \frac{1}{\| \psi\|^{p-2}} \left[
\frac{2(m-p) \| \psi\|_m^m}{m(p-2) \beta \|\psi\|^2}
\right]^{\frac{p-2}{m-2}}
$$
where $\sigma = \frac{{\mathcal S^{p/2}}|\Omega|^{1-p/2^{*}}}{p}$.
Since the right hand side in the last
inequality tends  as $m$ goes to $+\infty$ to 
$$
R \equiv \frac{\alpha}{4\sigma} \left[\frac{\|\psi
\|_\infty}{\|\psi \|}\right]^{p-2},
$$
it follows that  for every $\lambda >R$, there exists
$m_2(\lambda)>2^*$ such that if $m\geq m_2(\lambda)$, then
$T_m\|\psi\|_{X_{m}} > r_\lambda$, 
consequently,  every path $\gamma \in
\Gamma_m$ crosses the set $\|v\|_{X_{m}}=r_\lambda$, 
and we obtain
$$
c_m\geq \rho_{m,\lambda},
$$
(where $\rho_{m,\lambda}$ is defined in \eqref{weromin}). Therefore, the
geometrical hypotheses of the Mountain Pass Theorem are fulfilled
provided that
\begin{equation}\label{eq:ipolambda}
m \geq m_0\equiv \max \{ m_1(\lambda), m_2(\lambda)\} , \ \ \lambda
\geq \underline{\lambda}_{2} \equiv \max
\{\underline{\lambda}_{1} ,R \}.
\end{equation}
It is only left to show that  $J_m$ satisfies the Palais-Smale condition,
this is a straightforward argument; indeed, take a sequence 
$w_n $ satisfying
\begin{equation}\lbl{pscond}
\Jm(w_n)\stackrel{n\to +\infty}{\longrightarrow} c_m ,\qquad
J'_m(w_n)\stackrel{n\to +\infty}{\longrightarrow} 0 \quad \text{in
$X_m^*$},
\end{equation}
where $X_m^*$ denotes the dual space of $X_m$.
 As $\Jm$ is coer\-cive in $X_m$, the first convergence in
\eqref{pscond} implies that $w_n$ is bounded in $X_m$, then, up to a
subsequence, there exists $w\in X_m$ such that
\begin{equation*}
w_n \rightharpoonup  w\quad \text{ in $X_m$},
\end{equation*}
which means that $w_n \rightharpoonup  w$ in $L^m(\Omega )$ and in
$H_0^1(\Omega)$. Thus, by the Rellich theorem, $w_n$ is strongly convergent to
$w$ in $L^2(\Omega)$ and then, by interpolation, we deduce also the strong
convergence of $w_n$ in $L^r(\Omega)$ for all $r\in [2,m)$. Now, taking
$w_n-w\in H_0^1(\Omega) \cap L^{m}(\Omega ) \subset H_0^1(\Omega) \cap
L^{m'}(\Omega ) \subset X_m^*$ (since $m>2$) as test function in
\eqref{pscond}, we obtain
\[\begin{split}
\into M(x)\nabla w_n \cdot\nabla (w_n-w) +  \into |w_n|^{m-2} w_n (w_n - w)
\\-
 \lambda \into |w_n|^{p-2} w_n (w_n-w) \to 0 .
\end{split}\]
Observing that the strong convergence  of $w_n$ in $L^p(\Omega)$ implies that
the third adding term is tending to zero and, by subtracting the term 
$$
\into
M(x)\nabla w \cdot\nabla (w_n-w) +  \into |w|^{m-2} w (w_n - w),
$$ which is converging
to zero thanks to the weak convergence of $w_n$ in $\hsob$, we get
$$
\lim_{n\to +\infty}\into | \nabla (w_n-w)|^2
    +
 \into \left[ |w_n|^{m-2} w_n  -|w|^{m-2} w \right] (w_n - w)
        = 0 ,
$$
from which the strong convergence of $w_n$ to $w$ in $X_m$ is
deduced. Then, the Ambrosetti-Rabinowitz
theorem implies the existence of a mountain pass critical point
$z_m$ with
 critical level  $J_m(\zm)=c_m\geq 0>J_m(\um)$.
   Then, $\zm$ and $\um$ are distinct
solutions of problem \eqref{Prm}.
Finally, arguing as at  the end of Theorem \ref{esimini} it is possibile to obtain
that$z_{m}\geq 0$.
 \end{proof}
\begin{remark}
Let us observe that the result holds for every  $2<p\leq 2^{*}$ and for 
$m$ sufficiently large.
This assumption on $p$ is crucial to show that zero is a strict local minimum
so that $c_{m}>0$ for every $m$ fixed.
\end{remark}
\begin{remark}
Notice that  for $m$ sufficiently large $\rho_{m,\lambda}$ defined in \eqref{weromin} converges to zero as $m\to+\infty$, so that we do not 
have an immediate a priori bound from below on the action level of the
critical point $z_{m}$.
\end{remark}
\section{A priori Bounds}\label{bounds}
Let first show the following a priori  $L^{\infty}$-estimate that will be fundamental
in the asymptotical analysis.
\begin{prop}\lbl{stimenuova}

Assume conditions \eqref{alfa} and \eqref{pmu}. Then every positive solution
$w$ of Problem \eqref{Prm} satisfies the following estimate
\begin{equation}\label{stima:infty}
\| w\|_\infty \leq  \lambda  ^{\frac{1}{m-p}}.
\end{equation}
\end{prop}

 \begin{proof}
For every $t, \varepsilon >0$, let $\psi_\varepsilon$ the function given by
$$
 \psi_\varepsilon (s)  =
\left\{
\begin{array}{ll}
 0, & s<t,
             \\
             \\
             \frac{s-t}{\varepsilon} , & t<s<t+\varepsilon ,
             \\
             \\
             1, &t+\varepsilon \leq s.
\end{array}
 \right.
$$
Let us introduce the notation $E_t=\{ x\in \Omega :   w(x)>t\}$ and take $v=\psi_\varepsilon (w)$ as test function in
 $(\ref{equat})$ yielding
\[\begin{split}
   \int\limits_{E_t} w^{m-1} \psi_\varepsilon (w)
 \leq &
\int\limits_{\{x\in \Omega : t<w(x)<t+\varepsilon\}} M(x) \nabla w \cdot\nabla \psi_\varepsilon (w)
\\
&+ \int\limits_{E_t} w^{m-1} \psi_\varepsilon (w)
= 
 \lambda \int\limits_{E_t} w^{p-1} \psi_\varepsilon (w)
\leq
 \lambda \int\limits_{E_t} w^{p-1} .
\end{split}\]
Passing to the liminf as $\varepsilon $ goes to zero and applying
the Fatou Lemma and H\"older inequality,
\[\begin{split}
\int\limits_{E_t} w^{m-1}
 \leq
\lambda  \int\limits_{E_t} w^{p-1}
\leq
\lambda  \left[ \int\limits_{E_t} w^{m-1} \right]^{\frac{p-1}{m-1}} 
|E_t|^{\frac{m-p}{m-1}} .
\end{split}\]
That implies
\begin{equation}\label{wt}
 t^{m-1} |E_t|\leq \int\limits_{E_t} w^{m-1}
  \leq
 \lambda  ^{\frac{m-1}{m-p}}         |E_t| .
\end{equation}
Consequently, taking into account that $|E_t|\not= 0$ for every
$t\in (0,\|w\|_\infty)$ we obtain
$$
t\leq \lambda ^{\frac{1}{m-p}}, \ \ \forall
t\in (0,\|w\|_\infty)
$$
and hence \eqref{stima:infty} is proved.
 \end{proof}
\begin{remark}
The same conclusion of Proposition \ref{stimenuova} would hold
for every solution of Problem \eqref{Prm}; the hypothesis on positiveness is actually not needed. 
In addition,  notice that the nonlinearity $u^{m-1}$ satisfies all the  hypotheses of 
Theorem 1 in \cite{vaz}, then it actually holds that $u_{m}>0$ and $z_{m}>0$.
\end{remark}
\begin{remark}
The a-priori bound \eqref{stima:infty} prevents any blow-up phenomenon on
a sequence of positive solutions of Problem \eqref{Prm}, as we will see in the following.
\end{remark}
In addition  to Proposition \ref{stimenuova} we can also show the following.
\begin{prop}\lbl{stime}

Assume conditions \eqref{alfa} and  \eqref{pmu}. Then every solution
$w$ of Problem \eqref{Prm} satisfies the following estimates
\begin{align}
\lbl{stimam}
\|w\|_m^m &\leq \lambda ^{\frac{m}{m-p}}|\Omega|,
\\
\lbl{stimah1}
\alpha\|w\|^2+\|w\|_m^m &\leq |\Omega|
\lambda^{\frac{m}{m-p}}.
\end{align}
\end{prop}

 \begin{proof}
Taking $w$ as test function in \eqref{equat} and using \eqref{alfa}
we obtain
\begin{equation}\lbl{dis1}
\io|w|^m \leq \lambda  \io |w|^p
\leq \lambda 
 \left[ \io |w|^m \right]^{p/m}\!\!\!\!|\Omega |^{1-\frac{p}{m}},
\end{equation}
which implies that
$$
\left[ \io|w|^m \right]^{1-p/m} \leq \lambda  |\Omega
|^{1-\frac{p}{m}},
$$
and \eqref{stimam} is clearly deduced.
  Now, we choose again $v=w$ in
\eqref{equat}, to get from  H\"older inequality and  \eqref{stimam}
that
\begin{align*}
\alpha\io |\nabla w|^2+\io|w|^m &\leq
\lambda
  \left[ \lambda^{\frac{m}{m-p}}  \right]^{\frac{p}{m}}
  |\Omega |^{\frac{p}{m}} |\Omega |^{1-\frac{p}{m}}
 =\lambda^\frac{m}{m-p} |\Omega |
\end{align*}
 \end{proof}

\section{Asymptotical Analysis }\label{limit}
Let us start this section studying  the convergence of  the sequence of solutions
of minimum points $\{u_{m}\}$ 
\begin{teo}\label{limmini}
Assume \eqref{alfa} and \eqref{pmu}. There exists $\underline{\lambda}_{1}$ such that for every $\lambda>\underline{\lambda}_{1}$,
There exists  $u\in H^{1}_{0}(\Omega)\cap L^{\infty}(\Omega)$ such that
$u\not\equiv 0$, $u_{m}\rightharpoonup u$ weakly in $\hsob$, strongly in
every Lebesgue space and $u$ satisfies
\begin{equation}
\label{k}
u\in {\mathcal K}:=\{v\in H^{1}_0(\Omega):0\leq v(x) \leq 1 \},
\end{equation}
\begin{equation}\label{dv}
 \int_\Omega M(x)\nabla u    \cdot\nabla (v-u)    
 \geq \lambda \io \,u^{p-1}(  v - u ),\quad
 \forall\;v\in {\mathcal K} .
\end{equation}
In addition,  there exists
$g_{u}\in L^{\infty}(\Omega)$, such that
\begin{equation}\label{eq:gprop}
0\leq g_{u}\leq \lambda,\quad  g_{u}\not \equiv 0,\quad
g_{u}(x)[1-u(x)]=0, \;\text{  a. e. in $\Omega$},
\end{equation}
and it results
\begin{equation}\label{dv2}
 \int_\Omega M(x)\nabla u    \cdot\nabla \vfi+\into g_{u}\vfi=
 \lambda \io \,u^{p-1}\vfi ,\; \forall \vfi\in \hsob.
\end{equation}
\end{teo}

\begin{remark}
The properties of the  function $g_{u}$ expressed in \eqref{eq:gprop} show that
$g_{u}$ weights the set where $u\equiv1$, as $g_{u}(x)=0$ for almost
every $x\in \{x\in \Omega : u(x)<1\}$ and $g_{u}$ is not trivial.
\end{remark}

 \begin{proof}
We first apply Theorem \ref{esimini} to obtain a sequence of minimum
points $u_{m}$ of $J_{m}$ satisfying \eqref{stimamin}.
From \eqref{stimamin} and \eqref{stima:infty} we deduce that
there exists $u\in \hsob\cap \elle{\infty}$ with $0\leq u(x)\leq 1$
such that $u_{m}\rightharpoonup u$ weakly in $\hsob$, strongly in 
$\elle{q}$ for $q\in [1,+\infty)$. Moreover,  From
\eqref{stimamin} we also deduce that
\[
J_{\infty}(u_{m})=\frac12\into M(x)\nabla u_{m}\cdot \nabla u_{m}-
\frac{\lambda}p\into u_{m}^{p}\leq -\delta<0.
\]
Then, as $u_{m}\rightharpoonup u$ weakly in $\hsob$, and strongly in 
$\elle{q}$ for every $q\in [1,+\infty)$.
\[
0>-\delta\geq \liminf_{m\to+\infty}J_{m}(u_{m})\geq \frac12\into M(x)\nabla u 
\cdot \nabla u-\frac{\lambda}p\into u^{p}
\]
which implies that $u\not\equiv0$.
Let now $v$ be any element in ${\mathcal K} $, and let $\theta$ be any real number such
that $0 < \theta < 1.$ 
Using $\theta\,v-\um$ as test function in \eqref{equat}  we obtain
 \begin{equation}\label{2lug}
 \begin{split}
\int_\Omega M(x)\nabla \um  \cdot (\theta \nabla v - \nabla \um ) 
 &+ \int_\Omega   \,\um ^{m-1}(\theta v - \um )  
 \\
& = \lambda\io  \,\um^{p-1}(\theta v - \um ).    
\end{split}
\end{equation}
We write the second term   as
\[
\begin{split}
 \int_\Omega \um^{m-1}   (\theta v - \um ) =&
\int_{\{x  : 0\leq\um (x)  < \theta v\}} 
\um^{m-1}   (\theta v - \um )
\\
&+ 
\int_{\{x  : \um(x)   \geq \theta v\}} 
\um^{m-1}   (\theta v - \um ),
\end{split}    
\]
where the first term of the right hand side   is bounded by
$$
\um^{m-1}   (\theta v - \um )\leq
[\theta v]^{m-1}   (\theta v + \theta v )
\leq
2 \theta^{m}   ,
$$
  which tends to zero when $n$ tends to infinity because
$\theta < 1$; while
$$
 \um^{m-1}    (\theta v - \um ) \leq 0 \quad \hbox{on} \quad \{x :  \um (x) \geq\theta v\};
$$ 
and this implies that
$$
 {\limsup_{m\to \infty}} \int_\Omega  
\um^{m-1}   (\theta v - \um ) \leq 0.
$$
Now we write \eqref{2lug} as
\[
\begin{split}
\theta\int_\Omega M(x)\nabla \um \cdot   \nabla v   
 + 
  \int_{\{x : 0\leq\um(x)   < \theta v\}} 
\um^{m-1}   (\theta v - \um )
\\
 \geq \lambda \io \,\um^{p-1}(\theta v - \um )
 +\int_\Omega M(x)\nabla \um   \cdot \nabla \um , 
\end{split}\]
we pass to the limit as $n\to\infty$ and we 
observe that thanks to \eqref{alfa} we can exploit the weak 
lower semmicontinuity of the norm to obtain 
$$
\theta\int_\Omega M(x)\nabla u \cdot  \nabla v    
 \geq \lambda\io \,u^{p-1}(\theta v - u )
 +\int_\Omega M(x)\nabla u \cdot  \nabla u , 
$$
for any $v$ in ${\mathcal K} $ and any $\theta$ with $0 < \theta < 1$. Letting $\theta$
tend to 1 we get \eqref{dv}.

In order to prove the second part of the result, 
we take into account \eqref{stima:infty} and we deduce that
 there exists $g_{u}\in L^{\infty}(\Omega)$ such that,  (up to a subsequence),
 $\{(u_{m})^{m-1}\}\stackrel{*}{\rightharpoonup} g_{u}$ weakly-star  in $L^{\infty}(\Omega)$. As a consequence, we obtain that $g_{u}\geq 0$; in addition,
considering $\chi_{E}$ the characteristic function of
the set $E:=\{x\in \Omega : g_{u}>\lambda\}$  and exploiting  again \eqref{stima:infty} we get
\[
\lambda^{\frac{m-1}{m-p}}|E|\geq \into u_{m}^{m-1}\chi_{E}\quad \Longrightarrow\quad
\lambda |E|\geq \int_{E} g_{u}>\lambda |E|
\]
showing that 
\[g_{u}\leq \lambda.
\]
Taking $\vfi \in \hsob$ as  test function in 
\eqref{equat}and passing to the limit we get  that $u$ satisfies \eqref{dv2}.
Finally,
let us take 
as test function in \eqref{dv2} $v-u$ with $v\in {\mathcal K}$.
We obtain, that the equation
\[
\into M(x)\nabla u\cdot \nabla (v-u)+\into g_{u}(v-u)-\lambda\into
u^{p-1}(v-u)=0
\]
is satisfied for every $v\in {\mathcal K}$.
Then, using \eqref{dv} we deduce that 
\[
\into g_{u}(u-v)\geq 0\qquad \forall   v\in {\mathcal K} .
\]
Then, we can take a sequence $v_{j}\in {\mathcal K}$ such that
$v_{j}\to 1$ in $\elle{1}$ obtaining 
\[
\into g_{u}(u-1)\geq 0
\]
and this immediately implies that $g_{u}(1-u)\equiv 0$, as $g_{u}\geq0$ and $u\leq 1$.
In order to conclude, it is only left to show that $g_{u}\not \equiv 0$.
To this aim, we take into account that $J_{\infty}(u)\leq -\delta<0$, 
to get
\[
\lambda \into u^{p}\geq p\delta+\frac{p}2\into M(x)\nabla u\cdot \nabla u.
\]
On the other hand, choosing $\vfi=u$ in \eqref{dv2} we get
\[
\lambda \into u^{p}=\into g_{u}u+\into M(x)\nabla u\cdot \nabla u
\]
so that
\[
\into g_{u}u\geq p\delta+\left(\frac{p}2-1\right)\into M(x)\nabla u\cdot \nabla u
\]
and since $g_{u},\,u\geq 0,\, u\not\equiv 0$, this shows that $g_{u}\not\equiv 0$, or
equivalently $|\{x\in \Omega : u(x)=1\}|>0$.

 \end{proof}
\begin{remark}
Let us point out that the previous result holds for every $p>2$, so that the
nontrivial limit  solution $u$ exists even for $p>2^{*}$. 
\end{remark}
The previous results shows that, 
a similar phenomenon to the one observed in \cite{dalors}, \cite{bocmur} 
also occurs for this nonlinear problem.
Now, let us move to the study of the asymptotic behavior of the sequence of  the
critical points $z_{m}$, showing the following result.
\begin{teo}\label{limcrit}
Assume \eqref{alfa}, \eqref{pmz}. There exists $\underline{\lambda}\geq \underline{\lambda}_{1}$ such that for every $\lambda>\underline{\lambda}$,
there exists  $z\in H^{1}_{0}(\Omega)\cap L^{\infty}(\Omega)$ such that
$z\not\equiv 0$, $z_{m}\to z$ strongly in $\hsob$ and  in
every Lebesgue space. The function $z$ satisfies
\begin{equation}
\label{kz}
z\in {\mathcal K}:=\{v\in H^{1}_0(\Omega):0\leq v(x) \leq 1 \},
\end{equation}
\begin{equation}\label{dvz}
 \int_\Omega M(x)\nabla z    \cdot\nabla (v-z)    
 \geq \lambda \io \,z^{p-1}(  v - z ),\quad
 \forall\;v\in {\mathcal K} .
\end{equation}
In addition,  there exists
$g_{z}\in L^{\infty}(\Omega)$, such that
\begin{equation}\label{eq:gzprop}
0\leq g_{z}\leq \lambda, \qquad
g_{z}(x)[1-z(x)]=0, \;\text{  a. e. in $\Omega$},
\end{equation}
and it results
\begin{equation}\label{dvz2}
 \int_\Omega M(x)\nabla z   \cdot\nabla \vfi+\into g_{z}\vfi=
 \lambda \io \,z^{p-1}\vfi ,\; \forall \vfi\in \hsob.
\end{equation}
\end{teo}

 \begin{proof}
We first apply Theorem \ref{esimp} to obtain the existence of a sequence 
of critical points $z_{m}$ of $J_{m}$ for $m$  sufficiently large; then we
follow the same argument as in the proof of Theorem \ref{limmini} 
getting the existence of a function $z\in {\mathcal K}$ such that
$z_{m}\rightharpoonup z$ weakly in $\hsob$ and strongly in every
Lebesgue space; in addition $z$ satisfies \eqref{dvz}. 
Now, taking $z_{m}-z$ as test function \eqref{equat} yields
\begin{equation}\label{eq:1}
\begin{split}
\into M(x)\nabla z_{m}\cdot(\nabla z_{m}-\nabla z)&+\into |z_{m}|^{m-1}(z_{m}-z)
\\
&-\lambda \into  |z_{m}|^{p-1}(z_{m}-z)=0.
\end{split}
\end{equation}
Let us observe that
\[
\left|\into |z_{m}|^{m-1}(z_{m}-z)\right|\leq \into |z_{m}|^{m-1}|z_{m}-z|
\\
\leq \lambda^{\frac{m-1}{m-p}}
\into  |z_{m}-z|\to 0
\]
and the same argument shows that the last term in \eqref{eq:1} goes to zero.
Using these information in \eqref{eq:1} we get
\[
\into M(x)\cdot \nabla z_{m}(\nabla z_{m}-\nabla z)\to 0
\]
yielding
\[\begin{split}
\alpha\into |\nabla z_{m}-\nabla z |^{2}\leq 
&\into  M(x)\cdot \nabla z_{m}(\nabla z_{m}-\nabla z)
-
\into M(x)\cdot\nabla z (\nabla z_{m}-\nabla z)
\\=&
o(1)+\into M(x)\cdot\nabla z_{m}(\nabla z_{m}-\nabla z)
\end{split}\]
so that $z_{m}\to z$  strongly in $\hsob$. In order to show that $z\not\equiv 0$,
let us consider again the functional $J_{\infty}$ introduced in \eqref{def:Jinfty}.
Notice that there exists $\overline{T}$ and $\psi\in \hsob\cap \elle{\infty}$ such that
\[
J_{m}(\overline{T}\psi)<0 \quad \text{and }\quad J_{\infty}(\overline{T}\psi)<0
\] 
and the corresponding set of paths
$$
\Gamma_\infty=\left\{\gamma:[0,1]\to \hsob,\,:\gamma \text{ is continuous
and } \gamma(0)=0, \,J_{\infty}(\gamma(1))<0  \right\}.
$$
Moreover, we define the mountain pass  value
$$
c_\infty:=\inf_{\Gamma_\infty}\max_{[0,1]}J_\infty(\gamma(t)).
$$
We claim that $\Gamma_{m}\subset
\Gamma_{\infty}$. 
Indeed, take $\gamma\in \Gamma_{m}$ then $\gamma$ is evidently 
continuous in  $\hsob$ as it is continuous in $\hsob \cap L^{m}(\Omega)$,
and $\gamma(0)=0$. Moreover,
as  $\gamma(1)=T_{m}\psi$ with $J_{m}(T_{m}\psi)<0$, 
it is sufficient to notice that $J_{\infty}(T_{m}\psi)<J_{m}(T_{m}\psi)<0$
to obtain that $\gamma\in \Gamma_{\infty}$.
As a consequence, we get
\[
\max_{[0,1]}J_{m}(\gamma(t))\geq \max_{[0,1]}J_{\infty}(\gamma(t)) \quad (\text{because $J_{\infty}(\gamma)\leq J_{m}(\gamma)$ for every $\gamma$}).
\]
Then
\begin{equation}\label{stima:cm}
c_{m}=\inf_{\Gamma_{m}}\max_{[0,1]}J_{m}(\gamma(t))
\geq c_{\infty}=\inf_{\Gamma_{\infty}}\max_{[0,1]}J_{\infty}(\gamma(t))\geq \rho_{\infty} \; 
\end{equation}
because $\Gamma_{m}\subset\Gamma_{\infty}$  and
where $\rho_{\infty}>0$ can be obtained as 
\[
J_{\infty}(u)\geq \frac12\|u\|^{2}-C_{0}\|u\|^{p}\geq \rho_{\infty}>0
\]
for $\|u\|$ sufficiently small.
Therefore, 
\[
J_{m}(z_{m})= c_{m}\geq \rho_{\infty}>0.
\]
Then, exploiting \eqref{stima:infty} and recalling that $z_{m}\to z$ strongly in
$\hsob$ and in every Lebesgue space, we can  pass to the limit  to obtain 
that
\[
J_{\infty}(z)\geq \rho_{\infty}>0
\]
yielding that $z\not \equiv 0$.
Finally, the existence of $g_{z}$ satisfying \eqref{eq:gzprop} and \eqref{dvz2}
can be proved as in the proof of Theorem \ref{limmini}.
 \end{proof}
\begin{remark} \label{strongconv}
Notice that the strong $\hsob$ convergence can be also proved in an analogous way for the sequence of  minimum points $u_{m}$. However, while in the proof
of Theorem \ref{limcrit} this strong convergence is crucial to show that
$z$ is not trivial, in Theorem \ref{limmini} the weak convergence in $\hsob$ is 
sufficient to obtain that $u\not \equiv 0$.
\end{remark}
\begin{remark}\label{pcriti}
Notice that Theorem \ref{limcrit} holds for $p\leq 2^{*}$. This marks another difference
between the asymptotic behavior of $u_{m}$ and $z_{m}$. Indeed, $p\leq 2^{*}$  is not needed to pass to the limit, but it is crucial to obtain that $z\not \equiv 0$.
\end{remark}
\begin{remark}\label{gugz}
Let us point out that $g_{u}\not \equiv 0$ for every $p>2$. While, we cannot
show that $g_{z}\not \equiv0$ for any $2<p\leq 2^{*}$, we can only 
produce  some partial result which aim to show that $g_{z}$ could be trivial or not
depending on $p$ and on the domain.
\end{remark}
\begin{proposition}\label{pcrit}
Assume  $M(x)=Id$.
The following conclusions hold
\begin{enumerate}
\item
Let $\Omega$ be star-shaped, $p=2^{*}$, then $g_{z}\not \equiv 0$.
\item
Suppose that $\Omega=B_{1}(0)$, $N\geq 3$ and $p=2^{*}-\epsilon$. Then,
for $\epsilon >0 $ sufficiently small $g_{z}\not \equiv 0$.
\item
Suppose that $\Omega$ is a smooth bounded convex domain in  $\R^{2}$. Then,
for $p >2 $ sufficiently large $g_{z}\not \equiv 0$.
\end{enumerate}
\end{proposition}
\begin{proof}
{\bf Conclusion (1).}
Suppose by contradiction that $g_{z}\equiv 0$, then $z$ is a positive solution
to the problem
\[
\begin{cases}
-\Delta z=\lambda z^{2^{*}-1}& \text{in $\Omega$}
\\
z>0& \text{in $\Omega$}
\\
z=0& \text{on $\partial\Omega$}
\end{cases}\]
which cannot be by the Pohozaev identity (see \cite{poh} or \cite{struwe}).

{\bf Conclusion (2).}
It is known that (see \cite{adiyad, erbtan})  there exists a unique positive radially symmetric solution
$U_{\lambda}=\lambda^{-1/p-2}U$, where $U:B_{1}(0)\mapsto \R$ is
the unique positive solution of the Problem
\begin{equation}\label{eq:B1}
\begin{cases}
-\Delta U= U^{p-1} & \text{in }  B_{1}(0)
\\
\hskip0.62cm U=0 &\text{on } \partial B_{1}(0).
\end{cases}
\end{equation}
Moreover, as shown in \cite{atkpel} $\|U_{\lambda}\|_{\infty}\to +\infty$
as $\epsilon\to 0^{+}$. Then, for $\epsilon>0$ sufficiently small $U_{\lambda}>1$ is a set of positive measure, but $z\leq 1$ in the whole $\Omega$, so that
$g_{z}\not\equiv 0$.

{\bf Conclusion (3).}
We argue as in the proof of case  $(2)$: let $u_{p}$ be a positive solution of  
\begin{equation}\label{pLE}
\begin{cases}
-\Delta u_{p}= \lambda u_{p}^{p-1}& \text{in $\Omega$}
\\
u_{p}>0& \text{in $\Omega$}
\\
u_{p}=0& \text{on $\partial\Omega$}.
\end{cases}\end{equation}
with $\lambda=1$; then exploiting the result contained in \cite{demiangrospac}
$u_{p}$ is unique for $p$ sufficiently large. So that $u_{\lambda,p}=
\lambda^{-\frac1{p-2}}u_{p}$ is the unique, positive solution of Problem \eqref{pLE}.
But, as shown in \cite{adigro} 
\[
\|u_{\lambda,p}\|_{\infty}=\lambda^{-\frac1{p-2}}\|u_{p}\|_{\infty}\to \sqrt{e}\quad \text{as $p\to+\infty$}
\]
then, again $u_{\lambda,p}>1$ in a set of positive measure, for $p$ sufficiently
large, so that $g_{z}\not\equiv 0$ in this case too.

\end{proof}
Moreover, an analogous result holds for $p$ approaching 2 from above.
\begin{proposition}
Let $\Omega$ be any smooth bounded convex  domain in $\R^{N}$ with $N\geq 2$. Assume that  $M(x)=Id$. Then there exists $p_{0}>2$ such that
for every $p\in (2,p_{0})$ $g_{z}\not \equiv 0$.
\end{proposition}
\begin{proof}
The result can be obtained  following the argument in
Lemma 3.9 in \cite{pacdam}. Indeed, consider $u_{p}$ a positive 
solution of Problem \eqref{pLE} with $\lambda=1$, and 
$u_{\lambda,p}= \lambda^{-\frac1{p-2}}u_{p}$ solution of Problem \eqref{pLE}
for a  given $\lambda>\underline{\lambda}$.
Then, it is possible to show that, taking $p_{n}$ a sequence converging
to 2 from above, it results that $M_{n}:=\|u_{p_{n}}\|_{\infty}\to+\infty$
as $n\to +\infty$, and $M_{n}^{p_{n}-2}\to \lambda_{1}$; 
then 
\[
\lim_{n\to+\infty}\|u_{\lambda,p_{n}}\|^{p_{n}-2}_{\infty}=
\lim_{n\to+\infty}
\lambda^{-1}
\|u_{p_{n}}\|^{p_{n}-2}_{\infty}=\frac{\lambda_{1}}{\lambda}
\]
so that $\|u_{\lambda,p_{n}}\|_{\infty}\to +\infty$ too. Moreover, for $p\in (2,p_{0})$
this solution is unique, yielding again that
$g_{z}\not \equiv 0$, as $\|z\|_{\infty}\leq 1$.
 \end{proof}

Unfortunately, we  cannot show that $g_{z}\not\equiv 0$ for a
generic $p\in (2,2^{*})$, and we conjecture that this may depend on the domain
and on the exponent $p$, let us give some partial observations in this direction.

\begin{proposition}\label{gnontrivial}
Assume that  $M(x)=Id$. 
Moreover, suppose that  there exists $\lambda>0$ and $R>0$ such that
taking $\Omega=B_{R}(0)$,  the open ball of radius $R$ centred at zero,
the following conditions are satisfied
\begin{enumerate}
\item[(a)]
$\lambda>\underline{\lambda}$
\item[(b)]
$\left(\frac1{\lambda R^{2}}\right)^{\frac1{p-2}} U(0)>1$.
\end{enumerate}
Then $g_{z}\not \equiv 0$.
\end{proposition}
 \begin{proof}
Assume by contradiction that $g_{z}  \equiv 0$, then $z$ is a solution to the problem
\begin{equation}\label{plimz}
\begin{cases}
-\Delta z=\lambda z^{p-1} & \text{in }  \Omega
\\
z=0 &\text{on } \partial \Omega.
\end{cases}
\end{equation}
with $\Omega=B_{R}(0)$.
Then, taking into account  Theorem  (see \cite{adiyad, erbtan})
\[
z(x)=U_{\lambda,R}:=\left(\frac1{\lambda R^{2}}\right)^{\frac1{p-2}} U\left(\frac{x}R\right)
\]
where  $U$ is defined in \eqref{eq:B1}; indeed  $U_{\lambda,R}$
is defined in $B_{R}(0)$, vanishes at the boundary and solves
\[
-\Delta U_{\lambda,R}=\lambda^{-\frac1{p-2}}R^{-\frac2{p-2}-2}(-\Delta U)\left(\frac{x}R\right)
=\lambda^{-\frac1{p-2}}R^{-\frac2{p-2}-2}\left[U\left(\frac{x}R\right)\right]^{p-1}=\lambda 
U_{\lambda,R}^{p-1}
\]
In addition, recalling \eqref{kz} and hypothesis (b)
\[
1\geq z(0)=\max_{B_{R}(0)} z=\left(\frac1{\lambda R^{2}}\right)^{\frac1{p-2}} U(0)>1
\]
yielding the conclusion.
 \end{proof}
\begin{remark}\label{ending}
Unfortunately, we are  not able to give an example in which 
Proposition \ref{gnontrivial} applies showing that
$g_{z}\not \equiv 0$ in the subcritical regime as well. With this respect, let us observe that hypotheses 
(a) and   (b)  go in  opposite directions, as  (b) requires 
$\lambda$ sufficiently small, while in Theorem \ref{esimp} we have seen
if (a) is satisfied then
 $\lambda>\underline{\lambda}_{1}$, that
is
\[
\lambda>\underline{\lambda}_{1}:=\inf_{H^{1}_{0}(\Omega)\cap L^{\infty}(\Omega)}\Lambda(\vfi),\qquad
\Lambda(\vfi):=\frac{p}{2} \frac{\|\vfi \|^{2}\|\vfi\|^{p-2}_{\infty}}{\|\vfi\|_{p}^{p}}.
\]
For example,  
supposing that there exists $\vfi$ defined in $B_{1}(0)$ such that $\Lambda(\vfi)<\lambda$, then $\vfi_{R}(x):=\vfi(x/R)$ defined in $B_{R}(0)$ satisfies
\[
\Lambda(\vfi_{R})=\frac1{R^{2}}\frac{p}2 \frac{\|\vfi \|^{2}\|\vfi\|^{p-2}_{\infty}}{\|\vfi\|_{p}^{p}}=\frac1{R^{2}}\Lambda(\vfi)
\]
then, in order to have $g_{z}\not\equiv 0$ in $B_{R}(0)$, 
we would need
\[
\lambda>\frac1{R^{2}}\Lambda(\vfi),\quad\text{and}\quad
\left(\frac1{\lambda R^{2}}\right)^{\frac1{p-2}}U(0)>1
\]
which are  equivalent to
\[
\Lambda(\vfi)<\lambda R^{2}<\left[U(0)\right]^{p-2}.
\]
On  the other hand, let us observe that $\Lambda(U)=\frac{p}2\left[U(0)\right]^{p-2}$, so that choosing $\vfi=U$  would  imply that the assumptions of Proposition \ref{gnontrivial} cannot be satisfied.
So that finding an optimal $\vfi$  such that the interval $(\Lambda(\vfi), \left[U(0)\right]^{p-2})$ is not empty seems to be  a delicate question.
\end{remark}

\end{document}